\documentclass[a4paper,10pt]{article}
\newtheorem{thm}{Theorem}

\newtheorem{rem}[thm]{Remark}
\newtheorem{lem}[thm]{Lemma}
\newtheorem{prop}[thm]{Proposition}
\def \a{{\alpha}}

\def \d{{\delta}}
\def \e{{\varepsilon}}

\def \O{{\Omega}}

\font\phh=cmcsc10  at  8 pt

\def \A{{\cal A}}

\def \E{{\bf E}\, }

\def \N{{\bf N}}

\def \P{{\bf P}}
\def \q{{\quad}}

\def \qq{{\qquad}}
\def \R{{\bf R}}

\def \Z{{\bf Z}}

\def \noi{{\noindent}}
\def \L{{\Lambda}}

\def\qed{\hbox{\vrule height 6pt depth 0pt width
6pt}}
\def\cqfd{\hfill\penalty 500\kern 10pt\qed\medbreak}
 \font\phh=cmcsc10
at  8 pt \scrollmode \hfuzz =5 pt

\title{\bf A general strong law of large numbers
for additive arithmetic functions}
 \author{Istv\'an  Berkes and
   Michel Weber}

\begin{document}

\maketitle

  \begin{abstract}  Let $f(n)$ be a strongly  additive complex valued
  arithmetic   function. Under mild conditions on $f$,  we prove the
  following weighted strong   law of large numbers: if $ X,X_1,X_2,
  \ldots $ is any sequence of integrable i.i.d.\ random variables, then
  $$
  \lim_{N\to \infty}   { \sum_{n=1}^N   f(n)    X_n \over\sum_{n=1}^N   f(n)}
  \buildrel{a.s.}\over{=} \E X .
  $$
  \end{abstract}

\section{Introduction and main result}

Consider a strongly additive complex valued arithmetic function
$f(n)$, $n=1,2,\ldots$. Thus $f$ satisfies
\begin{eqnarray}f(mn) &=&f(m)+ f(n) \qq (m,n)=1,\cr
f(p^\a)&=&f(p), \qq p\  \hbox{\rm a prime},\ \ \a=1,2,\ldots
\end{eqnarray}
and it follows that
\begin{equation}  f(n)= \sum_{p|n}  f(p)   ,\label{(1.1)}
\end{equation}
so that $f$ is completely determined by its values taken over the
prime numbers.  We put
\begin{equation}  F(n) =\sum_{m\le n} f(m)\qq\quad
G(n) =\sum_{m\le n} |f (m)|^2. \label{(1.2)}
\end{equation}
Note that
\begin{equation}
F(n) =\sum_{p\le n} f(p)
\lfloor{n\over p}\rfloor, \q G(n) =\sum_{p\le n} |f
(p)|^2\lfloor{n\over p}\rfloor + 2\Re\Big\{\sum_{ 2\le  p< q \le n
}f (p)\overline{f (q)} \lfloor{n\over pq}\rfloor\Big\}.
\label{(1.3)}
\end{equation}
The general problem of  determining the order of magnitude of
additive arithmetic functions is a difficult task, and we refer to
the books Elliott \cite{E} and Kubilius \cite{ku} for a thorough
treatment.

\smallskip
In this work we are interested in the validity of the weighted
strong law of large numbers, when the weights are given by $f$.
More precisely, let ${\bf X}=\{X,X_m, m\ge 1\} $ be i.i.d.\ random
variables with basic  probability space $(\O,\A,\P)$  and such
that $\E |X|<\infty$.  We look for criteria for the weighted SLLN,
i.e.\ the relation
\begin{equation}\lim_{n\to \infty}   { \sum_{m=1}^n   f(m)  X_m
\over F(n)}\buildrel{a.s.} \over{=} \E X.  \label{(1.4)}
\end{equation}
 Such an SLLN  is a  delicate refinement
of the usual SLLN for i.i.d.\ random variables. Indeed, by
rewriting the sum in (\ref{(1.4)})  in the form
\begin{eqnarray} \sum_{m\le n}  f(m)X_m&=&\sum_{m\le n}  X_m
\sum_{p|m} f(p)=\sum_{p\le n} f(p) \sum_{p\le m\le n\atop p|m }X_m
\cr &=&\sum_{p\le n} f(p)\lfloor{n\over p}\rfloor \bigg({1\over
\lfloor{n\over p}\rfloor}\sum_{k=1}^{\lfloor{n\over p}\rfloor
}X_{kp} \bigg),
\end{eqnarray}
we see that  (\ref{(1.4)})  relies upon all  averages of the type
\begin{equation} \label{su}
{1\over m}\sum_{k=1}^{m}X_{kp}\qq\qq \hbox{$p$ prime   and $m\ge
2$},
\end{equation}
and in fact
it means that
$$\lim_{n\to \infty}   {1 \over \displaystyle{\sum_{p\le n} f(p)
\lfloor{n\over p}\rfloor }}\cdot\sum_{p\le n} f(p)\lfloor{n\over
p}\rfloor \bigg({1\over \lfloor{n\over
p}\rfloor}\sum_{k=1}^{\lfloor{n\over p}\rfloor }X_{kp} \bigg)\
\buildrel{a.s.}\over{=}\ \E X .
$$
Thus the validity of (\ref{(1.4)}) is intimately connected with
uniformity in the SLLN for the averages in (\ref{su}). Put
$$
A_n= \Big|\sum_{p\le n}{f(p)\over p}\Big|, \qq  B_n= \sum_{p\le n}
{|f (p)|^2 \over p}.
$$
If $f$ is real valued and nonnegative, then by  (\ref{(1.3)}) for
any $0<c\le 1/2$
\begin{equation}{c\over 2}nA_{c^2n}\le F(n)\le nA_n \qq {\rm and}\qq
G(n)\le n(B_n+  A_n^2).\label{(1.5)}
\end{equation}
So if $A_{2 n}\asymp A_n $, it follows that $F(n) \asymp nA_n$.
(Here, and in the sequel, $x_n\asymp y_n$ means
$0<\liminf_{n\to\infty}|x_n/y_n|\le \limsup_{n\to\infty}
|x_n/y_n|<\infty$.)

\smallskip\par  Naturally if $f$ is complex valued, the above
bound for $G(n)$  ceases to be true, because    the sum
 $$ \Re\Big\{\sum_{
  2\le  p< q \le n }f (p)\overline{f (q)} \lfloor{n\over pq}
  \rfloor\Big\},
  $$
is  in general no longer comparable to $A_n^2$.

 \medskip\par
In a recent work \cite{BW}, we studied the weighted SLLN when
$f\ge 0$ and proved the following result (see Theorem 1.1 in
\cite{BW}).
\begin{thm}\label{TheoremA}   Assume that $f\ge 0$ and
 \begin{equation} B_p\to\infty, \quad f(p)=o(B_p^{1/2})
\quad  {\textstyle as} \ p\to\infty. \label{(1.7)}
\end{equation}
  Then  (\ref{(1.4)}) holds.
\end{thm}

\medskip\par 
Condition  (\ref{(1.7)}) plays an important role in probabilistic
number theory as a nearly optimal sufficient condition for the
central limit theorem
\begin{equation} \label{ek}
\lim_{N\to \infty} \frac{1}{N} \# \{n\le N:
\frac{f(n)-A_N}{B_N^{1/2}} \le x\}= \frac{1}{(2\pi)^{1/2}}
\int_{-\infty}^x e^{-t^2/2} dt
\end{equation}
(see e.g.\ Elliott \cite{E}, Kubilius \cite{ku}.) Halberstam
\cite{Ha} proved that replacing the $o$ by $O$ in  (\ref{(1.7)})
the CLT (\ref{ek}) becomes generally false. Note that relation
(\ref{(1.7)}) implies the Lindeberg condition
\begin{equation} \label{lin}
\lim_{n\to \infty}{1\over B_n}\sum_{p<n\atop |f(p)|\ge \e
B_n^{1/2}}{f^2(p)\over p}= 0 \qq\hbox{for any $\e  >0$},
\end{equation}
and, under mild technical assumptions on $f$, condition
(\ref{lin}) is necessary and sufficient for the CLT (\ref{ek}),
see again Elliott \cite{E}, Kubilius \cite{ku}.

In \cite{BW} we also proved that (\ref{(1.7)}) implies the law of
the iterated logarithm corresponding to (\ref{(1.4)}) (see Theorem
1.2 in \cite{BW}). We futher indicated that if $f(p)$ does not
fluctuate too wildly, for instance if
\begin{equation}\label{fl} \sup_{n\le p,q\le n^2}{f(p)\over f(q)}=
{\cal O}(1),
\end{equation}
then Theorem \ref{TheoremA} remains valid under condition
(\ref{lin}).  We raised the question of the validity of Theorem
\ref{TheoremA} under the sole Lindeberg  condition. Recently,
Fukuyama and Komatsu \cite{FK} answered this question
affirmatively.
  \begin{thm}\label{TheoremB} Assume that $f\ge 0$ and the Lindeberg
 condition (\ref{lin}) is satisfied.
   Then  (\ref{(1.4)}) holds.
\end{thm}

Their approach is simple and elegant and is based on Abel
summation, and moreover it shows the interesting fact that the
Lindeberg condition implies
\begin{equation}\sum_{p>t} {f(p)\over p^2A_p}={\cal O}(1/t).\label{(1.8)}
\end{equation}
The estimates $F(n) \ge C_1 nA_n$, $G(n) \le C_2 nA_n^2$, which
are implied by  (\ref{lin}) (see for instance Lemma 2.1 in
\cite{BW}), are crucial in their proof, and their result remains
valid under these sole conditions. In particular, this is the case
if $A_{2 n}\asymp A_n $ and
\begin{equation}B_n^{1/2}={\cal O}(A_n ).\label{(1.9)}
\end{equation}
Actually, the condition
\begin{equation}F(n) \ge C_1 n\,\sup (A_n, B_n^{1/2})\label{(1.10)}
\end{equation}  would also suffice, since by   (\ref{(1.5)})  it implies
(\ref{(1.9)}) and thus $G(n) \le C_2 nA_n^2$. Condition
(\ref{(1.10)}) seems to be the relevant assumption in this
problem.  Note the interesting implication
\begin{equation} (\ref{(1.10)})\ \Longrightarrow \  (\ref{(1.8)}).
\label{(1.11)}
\end{equation}

 A typical example of application is the well-known von Mangoldt arithmetical
 function $\L$, which   is neither additive nor multiplicative.
Recall that $\L$ is  defined by
 \begin{equation}  \L(n)=
\cases{ \log p,  & \qq {\rm if}\quad $n=p^k$,\cr
 0,  & \qq {\rm otherwise}.}
\end{equation}
It is elementary to see that
$$
A_n= \sum_{p<n}{\L(p)\over p}\sim \log n, \qq  B_n= \sum_{p<n}
{\L^2(p) \over p}\sim \log^2 n .
$$
Since $A_n$ is slowly varying, on using  (\ref{(1.5)})  it follows
that  (\ref{(1.10)}) is valid. Therefore  (\ref{(1.4)}) holds when
the weights are given by von Mangoldt's function. Note also that
in this case (\ref{(1.8)}) reduces to the trivial estimate
 $\sum_{p>t}   p^{-2 }={\cal O}(1/t) $. By a result of Wierdl \cite{Wi},
 the result extends to the case when ${\bf X}$ is a   stationary
ergodic sequence in $L^p$, $p>1$ (and not if $p=1$, see
\cite{BM}), which is  a quite remarkable fact.
\medskip\par
It can   also be pointed out here that Abel summation alone
suffices to prove a result valid for general weights. In fact, the
proof in \cite{FK} yields the following

\begin{prop}\label{Proposition1.1}  Let
$f\ge 0$ be an arbitrary function and assume that
\begin{equation}\sum_{n> t}
 {G(n)\big|F (n+1)-F (n )\big| \over F^2(n)F (n+1)} = {\cal O}(1/t).
 \label{(1.12)}
\end{equation}
  Then  (\ref{(1.4)}) holds.
\end{prop}
\begin{rem}
{\bf (a)} \rm \ If there exists a   nondecreasing function $H(n)$
such that
\begin{equation}
F(n)\ge C_1 n H(n),   \q G(n) \le C_2 n H^2(n),
 \q
\big|F (n+1)-F (n )\big|  \le C_3H (n),  \label{(1.13)}
\end{equation}
   then condition    (\ref{(1.12)}) is satisfied and thus  (\ref{(1.4)}) holds.

\medskip \noi {\bf (b)} Under the Lindeberg condition we have
 $$\sum_{n> t}
 {G(n)\big|F (n+1)-F (n )\big| \over F^2(n)F (n+1)}\le C\Big({1\over t}
 +\sum_{p> t}
 { f(p ) \over p^2 A_p }\Big)= {\cal O}(1/t).  $$
\end{rem}
   \medskip\par\noi {\it Proof.} We have
 \begin{eqnarray} \#\big\{ n>t : F(n)\le tf(n)\big\} &\le&   t^2\sum_{n> t}
  {f^2(n)\over F^2(n)} \cr &= &  t^2\sum_{n> t}\Big(\sum_{t<k\le
n}f^2(k)\Big) \Big({1\over F^2(n)}-{1\over F^2(n+1)}\Big)
 \cr &\le &2 t^2\sum_{n> t}
 {G(n)\big|F (n+1)-F (n )\big| \over F^2(n)F (n+1)}
 \le   Ct.
\end{eqnarray}
Thus
$$\sup_{t>0}{1\over t}\#\big\{ n  : F(n)\le tf(n)\big\}<\infty,$$
which, by Lemma 2.1 (see Section 2) suffices to ensure
(\ref{(1.4)}). \cqfd

\bigskip\par
It is natural to ask about extensions of these results for complex
valued additive arithmetic functions.  As we mentioned after
(\ref{(1.5)}), the complex valued case requires a different
treatment. We do not know how to use Abel summation in this case.
Also, no estimate of type $F(n) \ge C_1 nA_n$, $G(n) \le C_2
nA_n^2$ or even $G(n) \le C_3n(B_n+  A_n^2)$ is available.
Further, the use of Abel summation
 leads to series  involving  $|f(n )|$, which are not related to
$$ A_{  n  }
  =\Big| \sum_{2\le p \le     n  }
{f(p)\over p}\Big| .$$ We will show, however, that a slight
modification in the use of the randomization argument introduced
in the proof of Theorem 1.1 in \cite{BW} allows in turn to prove a
rather general SLLN in this context.
\bigskip\par Let us first introduce some notation. Let    $\{\e_i, \,
i\ge 1\}$ denote a Bernoulli sequence defined on a probability
space $(\widetilde  \O,\widetilde \A,\widetilde \P)$,  with
partial sums $S_n=\e_1+\ldots+\e_n$. Let $\widetilde \E$ denote
the corresponding expectation symbol. Put
$$\eta= \sup \big\{ \rho>0\ :\ \widetilde   \P\big\{ \inf_{n\ge 1}{S_n\over n}
\ge\rho
 \big\}>0\big\}. $$
It is immediate to see that $\eta>0$. By the SLLN  we have ${S_n/
n} \to 1/2$ almost surely, so there is an integer $N \ge 3$ for
which
$$\widetilde \P\big\{ \inf_{n\ge N}{S_n\over n} \ge {1\over 3}
 \big\}\ge 2/3.
 $$
  Now
$$ {1\over 2}=\widetilde \P\{\e_1=1\}\le \widetilde \P\big\{ \inf_{n\le N}{S_n\over n} \ge
\inf_{n\le N}{\e_1\over n} \ge  {1\over N}
 \big\} . $$
Hence   $\widetilde \P\big\{ \inf_{n\ge 1}{S_n\over n} \ge {1\over
N}
 \big\}\ge 1/6$, which yields that $\eta\ge 1/N$.
\medskip\par
 
\medskip\par  Let $f$ be  a complex valued strongly additive arithmetic
function. We will prove the following result.
\begin{thm}\label{Theorem1.2}  Assume that there exists a nondecreasing
function $U:\N\to \R^+$ with $\lim_{n\to\infty} U(n)=\infty$ such
that
$$ c_1U(n)\le |F(n)|\ \le c_2 U(n)$$
for some positive constants $c_1, c_2$. Assume further that for
some $0<h<1/4$ we have
\begin{equation} \label{maincond}
\sup_{n^h < p \le n} |f(p)| \ll |F(\eta n)|/n, \quad A_{n^h} \ll
|F(\eta n)|/n,  \quad B_{n^h}^{1/2}  \ll  |F(\eta n)|/n.
\end{equation}
Then (\ref{(1.4)}) holds.
\end{thm}

\begin{rem} \rm  {\bf (a)}\q  Note that in condition (\ref{maincond}) we have
 $$ A_{  n^{h} }
  =\big| \sum_{2\le p \le     n^{h} }
{f(p)\over p}\big|\qq \hbox{instead of}\qq   \sum_{2\le p \le
n^{h} } \big|{f(p)\over p}\big|.$$ This was  made possible by
using a stronger estimate for divisors of Bernoulli sums than the
one used in \cite{BW}.

\medskip \noi {\bf (b)}\q If $f\ge 0$, one can take $U=F$. If
$F(n)\ge Cn\max(A_n, B_n^{1/2})$, then
 a sufficient condition for (\ref{maincond}) is
 $$\sup_{n^{h}<p \le n }  |f(p)|^2 \ll \max(A^2_{n^{h}}, B_{n^{h}}) $$
This is satisfied e.g.\ if $f$ is bounded.

\noi {\bf (c) }\ Condition (\ref{maincond}) can be replaced by a
slightly weaker condition of type (\ref{(1.12)}):
$$
\sum_{n\ge t} \frac{1}{|F(\eta n)|^2} \left(\sup_{n^h < p \le n}
|f(p)|^2 + A_{n^h}^2 + B_{n^h} \right) ={\cal O}(1/t).
$$
However, since the two conditions are close to each other and both
are probably far from being necessary and sufficient, it is
preferable to use the simpler assumption (\ref{maincond}).

\end{rem}

\section{Preliminaries}

In this section we formulate some lemmas needed for the proof of
Theorem  \ref{Theorem1.2}.

Let ${\bf X} =\{X_k, k\ge 1\}$ be i.i.d.\  random variables and
let ${\bf w}= \{w_k, k\ge 1\}$ be complex numbers with partial
sums $W_n=\sum_{k=1}^n w_k$, $n\ge 1$. We assume that
\begin{equation}
|W_n| \uparrow \infty, \qq n\to \infty.\label{(2.1)}
\end{equation}
Consider   the weighted averages
$$
M_n({\bf w}, {\bf X})={{1\over W_n}\sum_{k=1}^n w_k X_k }\qq
n=1,\ldots
$$
\begin{lem}\label{Lemma3}
We have $\lim_{n\to \infty} M_n({\bf w}, {\bf X})=0 $ almost
surely for every i.i.d.\ sequence $\bf X$  of nondegenerate,
centered, integrable random variables if and only if
$$
\limsup_{t\to\infty} {1\over t}\# \Big\{n:\ \left|{   W_n\over w_n
} \right|\le t\Big\} <\infty.
  $$
\end{lem}

Note that the last condition implies
$$\lim_{n\to\infty} \left|{   w_n\over W_n } \right|=0, $$
since, for any $\rho>0$, the number of integers $n$  such that $
|{   w_n/ W_n }  |>\rho$ is finite. The  characterization above is
due to   Jamison, Orey and Pruitt  (see Theorems 1 and 3 in
\cite{JOP}) under the additional fact that the weights $ w_k$ are
positive reals, in which case condition (\ref{(2.1)}) is trivially
satisfied.  As a matter of fact, the same proof allows  to work
with complex weights.
 
It would be natural to verify the conditions of Lemma \ref{Lemma3}
in order to prove Theorem \ref{Theorem1.2}, but for technical
reasons we were not able to do this. Instead, we will use the
following sufficient criterion for the weighted SLLN, also proved
(in the case of positive weights) in Jamison et al. \cite{JOP}.

\begin{lem}\label{lemmaL} Put
\begin{equation}  \qquad\qquad N(x)=
  \cases{
\#\{ k: \big|{W_k/ w_k}\big|\le x\}    \quad     &  \quad {\rm
if}\ $x\ge 1$, \  \  \  \  \  \, \cr 0            \quad    & \quad
{\rm if}      \  $0\le x<1$
 }    \label{(2.2)}
\end{equation}
and assume that $\E\,|X|<\infty$ and
\begin{equation}
 \E\,
    |X|^2  \Big(\int_{y\ge |X| }  { N(y ) \over y^3 }dy\Big)  <\infty.\label{(2.7)}
\end{equation}
Then we have $\lim_{n\to \infty} M_n({\bf w}, {\bf X})=0 $ a.s.
\end{lem}

\medskip
Again, the proof given in \cite{JOP} works in he complex case with
trivial changes. Let $\Psi$ denote the distribution function of
$X$ and let
$$
Y_k= X_k\!\cdot\!{\chi}  \big\{ |X_k| < |W_k/w_k| \big\}, \qquad
\zeta_k= |{ w_k /W_k}|\,\big(Y_k-\E Y_k\big).
$$
Following \cite{JOP}, we get
$$ \sum_{  k\ge 1} \P\{X_k\not =Y_k\}=\sum_{   k\ge 1}\int_{|v|\ge \big|{W_k\over
w_k}\big|} \Psi(dv) = \E N(|X|)
$$
and
$$\sum_{k=1}^\infty\E |\zeta_k|^2\le 4\int
    x^2  \Big(\int_{y\ge |x| }  { N(y ) \over y^3 }dy\Big)  \Psi(dx) .\label{(2.6)}
$$
As noted in \cite{JOP}, relation (\ref{(2.7)}) implies
$EN(|X|)<\infty$ and thus the lemma follows from the
Borel-Cantelli lemma and the Kolmogorov two series criterion.

Next we need a lemma on divisors of Bernoulli sums. Let
$d(n)=\#\{y :y|n\}$ be the divisor function.
 Consider the elliptic Theta function
\begin{equation}  \Theta (d,m)  =  \sum_{\ell\in \Z} e^{im\pi{\ell\over   d }
-{m\pi^2\ell^2\over 2 d^2}}.
     \label{(2.8)}
\end{equation}
The following lemma is   Theorem II from  \cite{We} which we
recall for convenience.
 \begin{lem}\label{Lemma4}      We have the following
uniform estimate:
\begin{equation}\sup_{2\le d\le n}\Big|\widetilde \P\big\{d|S_n\big\}- {\Theta(d,n)\over d}
 \Big|= {\cal O}\big((\log n)^{5/2}n^{-3/2}\big).
\label{(2.9)}
\end{equation}
 And
\begin{equation}
  \big|\widetilde \P\big\{d|S_n\big\}- {1\over d}  \big|\le  \cases{  C\Big((\log n)^{5/2}n^{-3/2}+
  { 1\over d}
 e^{ - {n \pi^2 \over 2d^2}} \Big)   &  \qq if $d\le \sqrt n$,\cr&\cr
      {C \over\sqrt n}   &  \qq if $ \sqrt n\le d\le n$,
}   \label{(2.10)}
\end{equation}
Further, for any $\a>0$
\begin{equation} \sup_{d<  \pi   \sqrt{   n \over 2\a\log n}}\big|\widetilde \P\big\{  d|   S_n
\big\}-{1\over d} \big|= {\cal O}_\e\big(n^{-\a+\e }\big) \quad
{\textstyle for \ all} \ \e>0.  \label{(2.11)}
\end{equation}
 and for any $0<\rho<1 $,
\begin{equation} \sup_{d<  (\pi/\sqrt 2) n^{(1-\rho)/2} }\big|\widetilde \P\big\{  d|   S_n \big\}
-{1\over d} \big|= {\cal O}_\e\big(e^{-(1-\e) n^\rho}\big),\q
{\textstyle for \ all} \ 0<\e<1.  \label{(2.12)}
\end{equation}
\end{lem}
\begin{rem} \rm By using the Poisson summation formula (see e.g.\
\cite{H}, p.\ 42)
\begin{equation}\sum_{\ell\in \Z} e^{-(\ell+\d)^2\pi x^{-1}}=x^{1/2} \sum_{\ell\in \Z}
e^{2i\pi \ell\d -\ell^2\pi x},\label{(2.13)}
\end{equation}
\par  \noi where $x$  is any real number and $0\le \d\le 1$,
  with the choices $x=\pi
n/(2d^2)$, $\d=n/(2d)$, we get
$${\Theta(d,n)\over d} ={1\over d}\sum_{\ell\in \Z} e^{i\pi n{\ell\over d} -n\pi^2
 {\ell^2\over 2d^2 }}= \sqrt{  { 2\over \pi
n}}\sum_{\ell\in \Z} e^{-2({n\over 2d}+\ell)^2 {d^2\over n}}
 .$$
Thus
 \begin{equation} \sup_{2\le d\le n}  \Big|\widetilde \P\big\{d|S_n\big\}-
 \sqrt{  { 2\over \pi n}}\sum_{\ell\in \Z} e^{-2( {n\over 2d} +\ell)^2
{d^2\over n}}  \Big|
  = {\cal O}\left({(\log n)^{5/2}\over n^{ 3/2}}\right).\label{(2.14)}
\end{equation}
Lemma \ref{Lemma4} was recently improved by the second named author for the range of values $d\ge \sqrt n$, one the basis of
these estimates.\end{rem}

  \section{Proof of Theorem \ref{Theorem1.2}}
We put
\begin{equation}L(t) =  \#  \{n:\   |F(n)| \le  t|f(n)| \}.\label{(3.1)}
\end{equation}
Since $\E|X|<\infty$, according to  Lemma \ref{lemmaL}, in order
to prove Theorem \ref{Theorem1.2}, it suffices to prove
\begin{equation}\label{(3.2)} \E\, |X|^2\int_{y\ge |X|}{L(y)\over y^3} dy<\infty.
\end{equation}
We use the same probabilistic trick as in \cite{BW}.
  We assume that the Bernoulli sequence $\{\e_i, \,
i\ge 1\}$ is defined on a  probability space   $  (\widetilde{\O},
\widetilde{\A}, \widetilde{\P})$,
  and denote by $ \widetilde{ \bf E}$   the corresponding  expectation symbol.
  Then, letting
  $F_\eta(n)=  \inf_{m\ge \eta n} |F(m)|$ we get
\begin{eqnarray}L(t)&\le &\#\big\{ n: |F(n)| \le t|f (n)|\big\}\le \#\big\{ n:
|F(S_n)|   \le t|f (S_n)|\big\}\cr &\le & \#\big\{ n: F_\eta(S_n)
\le t|f (S_n)|\big\}, \label{(3.3)}
\end{eqnarray}
 and this is true for  any $t>0$, simply because the graph of the
random walk $\{S_n, n\ge 1\}$ replicates all positive integers
with possible  multiplicities. If $\O_\eta=\{ S_n \ge \eta  n  \
{\rm for \ all} \ n\ge 1\} $ then   $\widetilde\P(\O_\eta)>0$.
Reading (\ref{(3.3)}) on $\O_\eta$ gives:
\begin{equation} L(t) \le \#\big\{ n : F_\eta(n) \le  t |f (S_n)|\big\} \quad \hbox{on
 $\O_\eta$  for all} \ t>0. \label{(3.4)}
\end{equation}
 But for all $t>0$
 \begin{eqnarray}
{1 \over t}  \#\big\{ n : F_\eta(  n)   \le t |f (S_n)|\big\}& \le
& 1+ {1 \over t}\#\big\{ n \ge t : F_\eta (n) \le t|f (S_n)|\big\}
\cr &= &1+ {1\over t} \sum_{n\ge t} \chi\{\, F_\eta^2(  n) \le
t^2|f  (S_n)|^2\}  \cr & \le & 1+t \sum_{n\ge t} {|f
(S_n)|^2\over F_\eta^2(  n) }  .   \label{(3.5)}
\end{eqnarray}
 We now prove the following lemma.
\begin{lem}\label{f(Sn)} There exists a constant $C_h$ depending on $h$ only such that
for any sufficiently large $n$
\end{lem}$$  \|f   (S_n) \|_{2,\widetilde \P}\le C_h\sup_{  n^{h}<p \le     n }
 |f(p)|+\big| \sum_{2\le p \le     n^{h} }
{f(p)\over p}\big|+C_\e\Big(\sum_{2\le p \le     n^{h} } {|f(p)|^2
\over p}\Big)^{1/2}.  $$

\noi {\it Proof.} We have
$$f   (S_n)  = \sum_{2\le p \le   S_ n }  f(p) \chi(
{p|S_n}) =  \sum_{2\le p \le     n }  f(p) \chi( {p|S_n}) ,
$$
and  given  any real $h$ with $0<h<1/4$
\begin{eqnarray}\Big| f   (S_n)-\sum_{2\le p \le     n^{h} }  f(p) \chi(
{p|S_n})\Big| &=&\big| \sum_{  n^{h}< p \le     n\atop p|S_n }
f(p) \big| \le \sum_{ n^{h}< p \le     n\atop p|S_n }  |f(p)|  \cr
&\le& C_h\sup_{  n^{h}<p \le     n }  |f(p)|. \label{(3.6)}
\end{eqnarray}
The last bound is justified by the fact that if $S_n$ admits $K$
different prime factors greater than $n^h$, then we have the
inequalities $n^{Kh}\le S_n\le n$; whence $Kh\le 1$.
 And so
\begin{equation}\label{depart}\Big| \|f   (S_n) \|_{2,\widetilde \P}-
\big\|\sum_{2\le p \le     n^{h} }  f(p) \chi(
{p|S_n})\big\|_{2,\widetilde \P}\Big|   \le C_h\sup_{  n^{h}<p \le
n }  |f(p)|.
\end{equation}
 Now denote by ${\bf 1}$ the function equal to $1$ everywhere on $\O$. Then
\begin{eqnarray}
& &  \Big|\big| \sum_{2\le p \le     n^{h} }  f(p)\widetilde
{\P}\{p|S_n\} \big|  -\big\| \sum_{2\le p \le     n^{h} }  f(p)
\chi( {p|S_n})\big\|_{2,\widetilde \P} \Big|\cr&= &
\Big|\big\|{\bf 1}\cdot\big( \sum_{2\le p \le     n^{h} }
f(p)\widetilde {\P}\{p|S_n\}
 \big)\big\|_{2,\widetilde
\P}  -\big\| \sum_{2\le p \le     n^{h} }  f(p) \chi(
{p|S_n})\big\|_{2,\widetilde \P} \Big|\cr &\le &   \Big\|
\sum_{2\le p \le     n^{h} }  f(p)\big( \chi(
{p|S_n})-\widetilde\P\{p|S_n\}\big)\Big\|_{2,\widetilde
\P}.\label{(3.7)}
\end{eqnarray}
Observe first by using Lemma \ref{Lemma4}
\begin{eqnarray*}
 \Big| \big| \sum_{2\le p \le     n^{h} }  f(p)\widetilde {\P}\{p|S_n\} \big|
 -\big| \sum_{2\le p \le     n^{h} }
{f(p)\over p}\big|  \Big|
  &\le & \Big|\sum_{2\le p \le     n^{h} }  f(p)\big(\widetilde {\P}\{p|S_n\}
  -{1\over
p}\big)  \Big| \cr &\le &  \sum_{2\le p \le     n^{h} } |
f(p)|\big |   \widetilde {\P}\{p|S_n\} -{1\over p}   \big | \cr
&\le & C_\e e^{-(1-\e) n^{1-2h}} \sum_{2\le p \le     n^{h} } |
f(p)| \
  .
 \end{eqnarray*}
Next by the Cauchy-Schwarz inequality
\begin{eqnarray}\label{cs}\Big( \sum_{2\le p \le     n^{h} }
  |f(p)|  \Big)^2&=&\Big( \sum_{2\le p \le     n^{h} }
 {|f(p)|\over \sqrt p}\cdot\sqrt p  \Big)^2 \le  \Big( \sum_{2\le p \le
    n^{h} }
 {|f(p)|^2\over  p}   \Big)\Big( \sum_{2\le p \le     n^{h} }
  p  \Big)\cr &\le & n^{2h}\Big( \sum_{2\le p \le     n^{h} }
 {|f(p)|^2\over  p}   \Big) .
\end{eqnarray}
Hence
\begin{eqnarray*}  e^{-(1-\e)
n^{1-2h}}  \sum_{ 2\le p \le     n^{h} } {|f(p)| } &\le &
e^{-(1-\e) n^{1-2h}}n^{  h}\Big( \sum_{ 2\le p \le     n^{h} }
{|f(p)|^2\over p }\Big)^{1/2}.
\end{eqnarray*}
Consequently,
\begin{eqnarray}  \label{(3.8)}
&& \Big| \big| \sum_{2\le p \le     n^{h} }  f(p)\widetilde
{\P}\{p|S_n\} \big|  -\big| \sum_{2\le p \le     n^{h} }
{f(p)\over p}\big|  \Big|
   \cr
  \cr &\le & C_\e e^{-(1-\e)
n^{1-2h}}n^{  h}  \Big( \sum_{ 2\le p \le     n^{h} }
{|f(p)|^2\over p }\Big)^{1/2} \cr &\ll & C_\e   \Big( \sum_{ 2\le
p \le     n^{h} } {|f(p)|^2\over p }\Big)^{1/2}.
 \end{eqnarray}
once $n$ is large enough, which we do assume. Thus (\ref{(3.7)})
and (\ref{(3.8)}) imply
\begin{eqnarray} \label{depart2}& &
\Big| \big\| \sum_{2\le p \le     n^{h} } f(p) \chi(
{p|S_n})\big\|_{2,\widetilde \P} -\big| \sum_{2\le p \le     n^{h}
} {f(p)\over p}\big|\Big|\cr &\le& \!\!\!  C_\e  \Big( \sum_{ 2\le
p \le     n^{h} } {|f(p)|^2\over p }\Big)^{1/2}\!\!\! +  \Big\|
\sum_{2\le p \le     n^{h} }  f(p)\big( \chi( {p|S_n})-{\widetilde
\P}\{p|S_n\}\big)\Big\|_{2,\widetilde \P}. \label{(3.9)}
\end{eqnarray}
 Clearly
   \begin{eqnarray}
  & & \widetilde{ \bf E}\,
\Big| \sum_{2\le p \le     n^{h} }  f(p) \big( \chi( {p|S_n})
-\widetilde\P\{p|S_n\}\big)\Big|^2 \cr  & =&\sum_{2 \le p  \le
n^{h}  }  |f (p ) |^2
 \big( \widetilde\P\{p|S_n\}-\widetilde\P\{p|S_n\}^2\big) \cr &\ & +
 2\Re\Big\{\sum_{2 \le
p<  q \le n^{h} }  f(p ) \overline{f(q )} \big(\widetilde{ \bf
P}\{  pq  |S_n  \}-\widetilde{ \bf P}\{   q |S_n  \}\widetilde{
\bf P}\{   q |S_n  \}\big)\Big\}
 .
\end{eqnarray}
But by Lemma \ref{Lemma4}
\begin{equation}\max\bigg\{  \big|\widetilde \P\big\{  p|   S_{ n} \big\}-{1\over p}
\big|, \big|\widetilde \P\big\{  pq|   S_{ n} \big\}-{1\over pq}
\big|\bigg\}\le  C_\e e^{-(1-\e) n^{1-4h}},
\end{equation}
whence
\begin{equation} \Big| \widetilde\P\{p|S_n\}\big(1- \widetilde\P\{p|S_n\} \big)-
{1\over p}(1-{1\over p})\Big| \le C_\e e^{-(1-\e)
n^{1-4h}},\label{(3.10)}
\end{equation}
and writing
\begin{eqnarray*} & &\widetilde{ \bf P}\{  pq  |S_n  \}
-\widetilde{ \bf P}\{   p |S_n  \}\widetilde{ \bf P}\{   q |S_n
\} \cr
 &  =&\big(\widetilde{ \bf P}\{  pq  |S_n  \}-{1\over pq}\big)
-\Big\{\big(\widetilde{ \bf P}\{   p |S_n  \}-{1\over p }\big)
\widetilde{ \bf P}\{   q |S_n  \}+{1\over p}\big(\widetilde{ \bf
P}\{   q |S_n  \}-{1\over q} \big) \Big\} ,
\end{eqnarray*}
we also get
\begin{equation} \Big|\widetilde{ \bf P}\{  pq  |S_n  \}-
\widetilde{ \bf P}\{   p |S_n  \}\widetilde{ \bf P}\{   q |S_n  \}
\Big|\le C_\e e^{-(1-\e) n^{1-4h}}. \label{(3.11)}
\end{equation}
By combining
 the previous relations it follows that
 \begin{eqnarray}& & \widetilde{ \bf E}\,
\Big| \sum_{2\le p \le     n^{h} }  f(p) \big( \chi( {p|S_n})
-\widetilde\P\{p|S_n\}\big)\Big|^2  \cr &\le &\sum_{2\le p \le
n^{h} } {|f(p)|^2 \over p}+ C_\e  e^{-(1-\e) n^{1-4h}} \Big(
\sum_{2\le p \le     n^{h} } {|f(p)| }\Big)^2 . \label{(3.12)}
\end{eqnarray}
By (\ref{cs})
 \begin{eqnarray*}  e^{-(1-\e)
n^{1-4h}} \Big( \sum_{ 2\le p \le     n^{h} } {|f(p)| }\Big)^2&\le
& e^{-(1-\e) n^{1-4h}}n^{ 2h}\Big( \sum_{ 2\le p \le     n^{h} }
{|f(p)|^2\over p }\Big)\cr &\ll & \sum_{ 2\le p \le     n^{h} }
{|f(p)|^2\over p }.
\end{eqnarray*}
Thus
\begin{equation} \widetilde{ \bf E}\,
\Big| \sum_{2\le p \le     n^{h} }  f(p) \big( \chi( {p|S_n})
-\widetilde\P\{p|S_n\}\big)\Big|^2  \ll \sum_{ n^{h}< p \le     n
} {|f(p)|^2 \over p}  . \label{(3.13)}
\end{equation}
By inserting (\ref{(3.13)}) into (\ref{(3.9)})  we arrive at
\begin{equation}
\Big| \big\| \sum_{2\le p \le     n^{h} } f(p) \chi(
{p|S_n})\big\|_{2,\widetilde \P} -\big| \sum_{2\le p \le     n^{h}
} {f(p)\over p}\big|\Big|  \le  C_\e\Big(\sum_{2\le p \le n^{h} }
{|f(p)|^2 \over p}\Big)^{1/2}
  .  \label{(3.14)}
\end{equation}
In view of (\ref{(3.6)}), the last relation implies
 \begin{equation}\label{d3}  \|f   (S_n) \|_{2,\widetilde \P}\le C_h
\sup_{  n^{h}<p \le     n }  |f(p)|+\big| \sum_{2\le p \le
n^{h} } {f(p)\over p}\big|+C_\e\Big(\sum_{2\le p \le     n^{h} }
{|f(p)|^2 \over p}\Big)^{1/2}
 .
\end{equation}
This completes the proof of Lemma \ref{f(Sn)}.\cqfd
\bigskip\par
We can now finish the proof of Theorem \ref{Theorem1.2}.  We get
from Lemma \ref{f(Sn)} and (\ref{(3.5)})
\begin{eqnarray*}& & \widetilde{ \bf E}\,  {1 \over t}
\#\big\{ n : F_\eta( n)
  \le t f (S_n)\big\} \le  1+ t \sum_{n\ge t} {\widetilde{\bf E}\,
|f   (S_n)|^2 \over  F^2_\eta( n)  }\cr & \le &1+C t \sum_{n\ge t}
{1\over F^2_\eta( n)}\bigg\{\sup_{  n^{h}<p \le     n }
  |f(p)|^2+\big| \sum_{2\le p \le     n^{h} }
{f(p)\over p}\big|^2  +  \sum_{2\le p \le     n^{h} }
{|f(p)|^2 \over p}  \bigg\}    .  
\end{eqnarray*}
 On using assumption  (\ref{maincond}), we deduce
$$
\sup_{t>0}\widetilde{ \bf E}\,  {1 \over t}
\#\big\{ n : F_\eta( n)   \le t |f (S_n)|\big\}\le
C. 
$$
It follows that
\begin{eqnarray*}
 \E \, \widetilde{ \bf E}\, X^2\int_{y\ge |X|}{ \#\big\{ n : F_\eta( n)
\le y|f(S_n)|\big\}\over y^3} dy &\le &C \E \, X^2\int_{y\ge |X|}{
1\over y^2} dy\cr &\le& C  \E \, |X|<\infty. \label{(3.18)}
\end{eqnarray*}
And  in view of (\ref{(3.3)}), (\ref{(3.5)}) and Fubini's theorem
\begin{eqnarray*}& &\widetilde{ \bf E}\,  \chi(\O_\eta) \cdot \E \,
X^2\int_{y\ge |X|}{L(y)\over y^3} dy\cr  &\le&\E \, \widetilde{
\bf E}\, X^2\int_{y\ge |X|}{ \#\big\{ n : F_\eta( n) \le y
|f(S_n)|\big\}\over y^3} dy \cr &\le& C  \E \, |X|<\infty.
 \label{(3.19)}
\end{eqnarray*}
Since
$$
\widetilde{\bf E}\, \chi(\O_\eta) \cdot \E \, X^2\int_{y\ge
|X|}{L(y)\over y^3} dy =\P\{\chi(\O_\eta)\}  \E \, X^2\int_{y\ge
|X|}{L(y)\over y^3} dy,
$$
relation (\ref{(3.2)}) follows, completing the proof.


{
 \bigskip\par\noi {\phh Istvan Berkes, Graz University of Technology,
Institute of Statistics, M\"unzgraben-strasse 11, A-8010 Graz,
Austria.
 E-mail:\ \tt    berkes@tugraz.at}
\medskip\par
\noi {\phh Michel Weber:  IRMA, Universit\'e Louis-Pasteur et
C.N.R.S.,   7  rue Ren\'e Descartes, 67084 Strasbourg Cedex,
France. \noi E-mail: \  \tt weber@math.u-strasbg.fr}

\end{document}